%-----------------------------------------------------------------------------
% generic.tex
%                    
% Beginn: Dezember 2012
%------------------------------------------------------------------------------
\magnification=\magstep1   
\input amstex
\UseAMSsymbols
\input pictex
%\input german 
%\hoffset=0truecm \voffset=0truecm 
\vsize=23truecm
\NoBlackBoxes
\parindent=20pt
  \font\rmk=cmr8 
\font\ttk=cmtt8
  
%\hrule height 2pt \vskip 3pt \hrule \bigskip\bigskip\bigskip

\def\bdim{\operatorname{\bold {dim} }}
\def\bDim{\operatorname{\bold {Dim} }}
\def\mod{\operatorname{mod}}
\def\Mod{\operatorname{Mod}}

\def\Hom{\operatorname{Hom}}
\def\End{\operatorname{End}}

\def\op{{\text{op}}}
\def\Ext{\operatorname{Ext}}
\def\rad{\operatorname{rad}}

    %Fuer Text innerhalb von $$...$$. Verwendung: \T{oder}
    %Fuer Text innerhalb von $$...$$. Verwendung: \T{oder}
   
          %Zum Einruecken, Verwendung: \E{dieses wird eingerueckt}
\def\arr#1#2{\arrow <1.5mm> [0.25,0.75] from #1 to #2}

\vglue2truecm

\centerline{\bf Generic representations of wild quivers.}
	\bigskip
\centerline{Claus Michael Ringel}
	\bigskip\bigskip
{\narrower\narrower
{\bf Abstract.} Let $\Delta$ be a wild connected directed quiver.  We show that any generic 
representation $M$ 
is the union of its subrepresentations of finite length which are regular. As
a consequence, we see that the direct limit closure of the preprojective component 
does not contain any 
generic module.
\par}

	\bigskip
Let $\Delta$ be a finite quiver which is connected and directed. 
We
consider representations of $\Delta$ with coefficients in the field $k$, or,
what is the same, $\Lambda$-modules, where $\Lambda = k\Delta$ is the path algebra of 
$\Delta$. Note that the path algebra of a quiver is hereditary, and since we assume
that $\Delta$ is directed, $\Lambda$ is finite-dimensional.

Given a module $M$, let $E(M) = \End(M)^\op$. We may consider $M$ as an $E(M)$-module.
The module $M$ is said to be {\it endo-finite,} provided it is of finite length when
considered as an $E(M)$-module. 
	\medskip 
{\bf Theorem.} {\it Let $\Delta$ be a connected directed finite quiver. Let
$M$ be an indecomposable endo-finite representation of infinite length. If $\Delta$ is
wild, then $M$
is the union of its subrepresentations of finite length which are regular.}
	\medskip 
{\bf Corollary.} {\it Let $\Delta$ be a connected directed finite quiver. If $\Delta$
is wild, then the direct limit closure of the preprojective component contains no
indecomposable endo-finite modules of infinite length.}
	\medskip 
This answers a question raised by Henning Krause, see [K] for consequences. 
Indecomposable endo-finite modules of infinite length are sometimes called {\it generic} modules. 
If $\Delta$ is representation-finite, then 
there is no generic module. If $\Delta$ is tame, then there is a unique
generic module $M$ and $\Hom(R,M)= 0$ for any finite length module $R$ which is regular.
In this case, $M$ is in the direct limit closure of the preprojective component.
	
The main tool for the proof is the endo-length vector $\bDim M$ of an endo-finite module $M$
which will be introduced in section 2. In section 3, we first will show that given a generic
module $M$, there is a regular module $R$ with $\Hom(R,M) \neq 0$, and we use this result
in order to provide the proof of the Theorem. 
	\medskip
{\bf Acknowledgment.} The author is indebted to O\. Kerner and H\. Krause for helpful
comments concerning the presentation of the results. 
	\bigskip\bigskip
\vfill\eject
%=========================================================================
{\bf 1. Preliminaries on finite length modules.}
	\medskip 
Let $\Delta_0$ be the set of vertices, $\Delta_1$ the set of arrows
of the quiver $\Delta$. Given an arrow $\alpha$, let $s(\alpha)$ be its starting vertex,
$t(\alpha)$ its terminal vertex. Let $\Lambda = k\Delta$ and 
$\Mod \Lambda$ be the category of $\Lambda$-modules,
and $\mod \Lambda$ the full subcategory of all $\Lambda$-modules of finite length.
We denote by $K_0(\Lambda)$ the Grothendieck group
of $\mod \Lambda$ (with respect to all exact sequences), we may identify
it with the set of functions $\Delta_0 \to \Bbb Z$, thus 
with $\Bbb Z^n$ if $\Delta_0$ has cardinality $n$.
Given a module $M$ of finite length,
the corresponding element in $K_0(\Lambda)$ is called its dimension vector and
denoted by $\bdim M$, the coefficient $(\bdim M)_i$ for $i\in \Delta_0$ is the
Jordan-H\"older multiplicity of the simple module $S(i)$ in $M$. 
For $x = (x_i)_i,\ y = (y_i)_i\in K_0(\Lambda)$, one defines
$$
 \langle x,y\rangle = \sum_{i\in \Delta_0} x_iy_i 
 - \sum_{\alpha \in \Delta_1} x_{s(\alpha)}y_{t(\alpha)},
$$
and one obtains in this way an integral bilinear form $\langle-,-\rangle$ on $K_0(\Lambda)$.
Since 
$$
 \langle \bdim X,\bdim Y\rangle = \dim_k \Hom(X,Y) - \dim_k\Ext^1(X,Y),
$$
(and $\Ext^i(X,Y) = 0$ for $i\ge 2$)
the form $\langle-,-\rangle$ is called the {\it Euler form}.

We denote by $\tau$ the Auslander-Reiten translation on the category of $\mod \Lambda$. 
Recall that a
indecomposable module $X$ of finite length 
is called preprojective or preinjective 
provided $\tau^tX = 0$ or $\tau^-tX$, respectively, for 
some natural number $t$. A module will be said to be {\it regular} 
 provided it has finite length
and has no indecomposable direct summand which is preprojective or preinjective. 
We denote by $\Phi$ the Coxeter transformation on $K_0(\Lambda)$, it is a
linear transformation and $\bdim \tau X = \Phi\bdim X$ for any 
indecomposable non-projective module $X$ of finite length.
	\medskip
We say that a module $X$ is {\it in general position} provided $X = X'\oplus X''$
implies that $\Ext^1(X',X'') = 0.$ 
If $x$ is a dimension vector, then 
the modules $X$ with $\bdim X = x$
such that $\dim_k\End(X)$ is minimal, are in general position. 
	\medskip 
{\bf Lemma 1.} {\it The following conditions are equivalent for $x\in K_0(\Lambda).$
\item{\rm (i)} $\Phi^t(x) \ge 0$ for all $t\in \Bbb Z$.
\item{\rm(ii)} The finite length modules in general position with dimension vector $x$
    are regular.
\item{\rm(iii)} There exists a regular module with dimension vector $x$.\par}
	\medskip
Note that in (iii) we cannot expect that there exists a regular module which is indecomposable,
typical examples are the elements in the $\Phi$-orbits of
$(1,1,0,1,1)$ and $(1,1,3,1,1)$ for the quiver
$$
{\beginpicture
\setcoordinatesystem units <1cm,1cm>
\multiput{$\circ$} at 0 0  1 0  2 0  3 0  4 0 /
\arr{0.8 0.1}{0.2 0.1}
\arr{0.8 -.1}{0.2 -.1}
\arr{1.8 0}{1.2 0}
\arr{2.8 0}{2.2 0}
\arr{3.8 0.1}{3.2 0.1}
\arr{3.8 -.1}{3.2 -.1}
\endpicture}
$$
Let us call $x\in K_0(\Lambda)$ {\it regular} provided the equivalent conditions of 
Lemma 1 are satisfied.
	\medskip
Proof of Lemma 1: (i) $\implies$ (ii): Let $X$ be a module in general position and 
$\bdim X = x$. Let us assume that $X$ has an non-zero preprojective direct
summand, say let $X = X'\oplus \tau^{-t}P$, where $P$ is a non-zero projective
module, $t\ge 1$. We can assume that $X'$ has no non-zero direct summand of the
form $\tau^{-s} P'$ with $P'$ projective and $0 \le s \le t.$ With $X$ also 
$\tau^t X = \tau^tX'\oplus P$ is in general position, thus
$$
 0 = \Ext^1(\tau^tX',P) = D\Hom(P,\tau^{t+1}X')
$$
(here, $D$ denotes the $k$-duality). 
This means that $(\bdim \tau^{t-1}X')_j = 0$ for every vertex $j$ such that the
corresponding indecomposable projective module
$P(i)$ is a direct summand of $P$. Now 
$$
 \Phi^{t+1}\bdim X = \Phi^{t+1}\bdim X'- \bdim \nu P = 
 \bdim \tau^{t-1}X' - \bdim \nu P,
$$ 
where $\nu$ is the Nakayama functor (it sends $P(i)$ to the
corresponding indecomposable injective module $I()$). Since $P$
is non-zero, there is a vertex $j$ such that $P(j)$ is a direct summand
of $P$, therefore $I(j)$ is a direct summand of $\nu P$. Since $(\bdim \tau^{t-1}X')_j = 0$
and  $(\bdim \nu P)_j > 0$, we see that $\Phi^{t+1}\bdim X$ has a negative coefficient.

Using duality, we similarly see that $X$ has no indecomposable preinjective direct 
summand. Therefore $X$ is regular.

(iii) $\implies$ (i) follows from the fact that for a regular module $R$, also
$\tau R$ and $\tau^-R$ are regular and $\Phi(\bdim R) = \bdim \tau R$
and $\Phi^{-1}(\bdim R) = \bdim \tau^- R$.
	\bigskip
{\bf Lemma 2.} {\it Let $x$ be a non-zero regular element of $K_0(\Lambda)$.
Then $\langle \Phi^{-t}x,x\rangle > 0$ for $t \gg 0.$}
	\medskip
Proof: Let $X$ be a non-zero regular module with $\bdim X = x$ and apply Baer [B],
Proposition 2.1 for $X = S$.
 	\bigskip\bigskip 
%=========================================================================
{\bf 2. The endo-length vector $\bDim M$ of an endo-finite module $M$.}
	\medskip 
In [R3], for any endo-finite modules $M$ whose endomorphism ring is 
a division ring an element $\bDim M$ in $K_0(\Lambda)$ has been defined
(in a similar setting, Lenzing [L] has called such an invariant 
$\bDim M$ the ``characteristic class'' of $M$). 
We extend the definition from [R3] to arbitrary endo-finite modules.
 Note that for any representation $M = (M_i,M_\alpha)$ of $\Delta$, the vector spaces $M_i$ are
$E(M)$-modules. If $M$ is endo-finite, all the
$E(M)$-modules $M_i$ have finite length $|{}_{E(M)}M_i|$ and 
$|{}_{E(M)}M| = \sum_{i\in \Delta_0} |{}_{E(M)}M_i|$.
We introduce the {\it endo-length vector} $\bDim M$ as the function 
$\bDim\:\Delta_0 \to \Bbb Z$ with
$$
 (\bDim M)_i = |{}_{E(M)}M_i|,
$$
thus $\bDim M$ is an element of $K_0(\Lambda).$ 
	\medskip 
Remark. If $M$ is an indecomposable module of finite length, then both $\bdim M$ and $\bDim M$
are defined, and are multiples of each other, namely we have 
$\bdim M = \dim_k\underline{\End} M
\cdot \bDim M$, where $\underline{\End} M = \End M/\rad\End M.$ Of course, in case $k$
is an algebraically closed field, $\underline{\End} M = k$ for all indecomposable modules $M$,
thus $\bdim M = \bDim M$. But if $k$ is not algebraically closed, then already for the
Kronecker quiver there exist indecomposable modules of finite length with 
$\dim_k\underline{\End} M > 1.$

	\medskip
Recall that Bernstein-Gelfand-Ponomarev [BGP] have defined reflection functors.
If $i$ is a vertex of $\Delta$, the quiver $\sigma_i\Delta$ is obtained from $\Delta$
by changing the orientation of all the arrows starting or ending in $i$.
Given a sink $i$ of the quiver, there is the reflection functor $\sigma_i\:\Mod k\Delta
\to \Mod k(\sigma_i\Delta)$, it provides an equivalence between the full subcategory of
$\Mod k\Delta$ of all $k\Delta$-modules without direct summands $S(i)$ and
the full subcategory of
$\Mod k(\sigma_i\Delta)$ of all $k(\sigma_i\Delta)$-modules without  direct summands $S(i)$.
Similarly, for $i$ a sink, there is the corresponding reflection functor 
$\sigma_i\:\Mod k\Delta \to \Mod k(\sigma_i\Delta)$. In addition, we also denote by
$\sigma_i$ the
reflection $\sigma_i\:\Bbb Z^n \to \Bbb Z^n$ such that $\bdim \sigma_iM = 
\sigma_i\bdim M$ for any finite length module $M$ without a direct summand of the form $S(i)$
(and $i$ a sink or a source). 
	\medskip
{\bf Lemma 3.} {\it Let $M$ be a generic module. If $i$ is a sink of the quiver $\Delta$,
then $\bDim \sigma_iM = \sigma_i\bDim M.$
If $i$ is a source of the quiver $\Delta$,
then $\bDim \sigma_i^-M = \sigma_i\bDim M.$}
	\medskip
Proof: We only discuss the case of $i$ being a sink. 
By definition, $(\sigma_i M)_i$ is the kernel of
the map 
$$
 \bigoplus_{t(\alpha)=i} M_{s(\alpha)} \to M_i
$$
whose  restriction to $M_{s(\alpha)}$ is $M_\alpha$. Since $M$ has no direct summand of
the form $S(i)$, this map is surjective. Also, this is an $E(M)$-module
homomorphism. Altogether, we see that we deal with the exact sequence
$$
 0 \to (\sigma_iM)_i \to \bigoplus_{t(\alpha)=i} M_{s(\alpha)} \to M_i \to 0
$$
of $E(M)$-modules. Looking at the length of these $E(M)$-modules, we have
$$
 |{}_{E(M)}(\sigma_i M)_i| = \sum _{t(\alpha)=i} |{}_{E(M)}M_{s(\alpha}| 
  - |{}_{E(M)}M_i| = (\sigma_i\bDim M)_i.
$$
This shows that $\bDim \sigma_i(M) = \sigma_i\bDim M.$
	\medskip 
If we label the vertices of $\Delta$ as $\Delta_0 = \{1,2,\dots,n\}$ so that there is
no arrow $i \to j$ for $i\le j$, then the composition $\Phi = \sigma_n\cdots\sigma_1$
of the reflection functors $\sigma_i$ is  defined and is called 
{\it Coxeter functor,} the corresponding composition $\Phi = \sigma_n\cdots\sigma_1$
of the reflections $\sigma_i\:\Bbb Z^n \to \Bbb Z^n$ 
is the Coxeter transformation mentioned already.
	\medskip 
{\bf Lemma 4.} {\it Let $M$ be a generic module. Then $\Phi M,\ \Phi^{-1}M$ are 
generic modules and }
$$
 \bDim\Phi M = \Phi\bDim M, \quad \bDim\Phi^{-1} M = \Phi^{-1}\bDim M.
$$ 
	\medskip
Proof. This follows immediately from the fact that 
$M$ has no non-zero projective direct summands.
	\medskip 
{\bf Lemma 5.} {\it  Let $M$ be a generic module. Then
$\bDim M$ is a non-zero regular element of $K_0(\Lambda).$}
	\medskip
Proof: Let $x = \bDim M$. According to Lemma 4, all the vectors $\Phi^tx$ are
endo-length vectors of non-zero modules, thus non-negative. 
	\bigskip  
{\bf Lemma 6.} {\it Let $X$ be of finite length and $M$ endo-finite. 
Then there is an exact sequence of $E(M)$-modules}
$$
 0 @>>> \Hom(X,M) \to 
 \bigoplus_i \Hom_k(X_i,M_i) @>\delta_{XM}>> \bigoplus_{\alpha}
 \Hom_k(X_{s(\alpha)},M_{t(\alpha)}) @>>> \Ext^1(X,M) @>>> 0.
$$
	\medskip
Proof. Let us refer to [R1] where the special case of both $X,M$ being of finite length
has been considered, but without taking into account the $E(M)$-module structure.
As in the special case, the map $\delta_{XM}$ is defined in general as follows: it sends
an element  $f = (f_i)_i$ with $f_i\in \Hom_k(X_i,M_i)$ to $\delta_{XM}(f)$ with
components $\left(\delta_{XM}(f)\right)_\alpha = 
f_{t(\alpha)}X_\alpha - M_\alpha f_{s(\alpha)}.$ It is clear that this map is an
$E(M)$-homomorphism. Now it is trivial to verify that the kernel of $\delta_{XM}$ 
is just $\Hom(X,M)$, thus only the assertion that 
the cokernel of $\delta_{XM}$ is equal to $\Ext^1(X,M)$ has to be shown. The
proof given in [R1] remains true in our more general setting.
	\medskip 
As an immediate consequence we obtain:
	\medskip 
{\bf Lemma 7.} 
{\it Let $X$ be of finite length and $M$ endo-finite. 
Then }
$$
 \langle \bdim X,\bDim M\rangle =  |{}_{E(M)}\Hom(X,M)| -  |{}_{E(M)}\Ext^1(X,M)|.
$$
	\medskip
Proof: Let $\bdim X = x$. Note that for any $d$-dimensional vector space $V$, we have
$|{}_{E(M)}\Hom_k(V,M_i)| = d |{}_{E(M)}M_i|.$
The direct sum decompositions 
$\bigoplus_i \Hom_k(X_i,M_i)$ and $\bigoplus_\alpha  \Hom_k(X_{s(\alpha)},M_{t(\alpha)})$
are direct sums of $E(M)$-modules, thus, using Lemma 6, we have
$$
\align
 |{}_{E(M)}\Hom(X,&M)| -  |{}_{E(M)}\Ext^1(X,M)| \cr
 &=  
 \sum_i|{}_{E(M)}\Hom_k(X_i,M_i)| -  \sum_{\alpha}|{}_{E(M)}\Hom_k(X_{s(\alpha)},
   M_{t(\alpha)})| \cr
 &=  
 \sum_ix_i|{}_{E(M)}M_i| -  \sum_{\alpha}x_{s(\alpha)}|{}_{E(M)}M_{t(\alpha)}| \cr
 &=\langle \bdim X,\bDim M\rangle 
\endalign
$$
	\medskip 
This implies:
	\medskip 
{\bf Lemma 8.} {\it Let $X$ be of finite length and $M$ endo-finite. 
If $\langle \bdim X,\bDim M\rangle > 0,$ then $\Hom(X,M) \neq 0$.}
	\medskip
We should note that for Lemma 8, we only need the trivial part of Lemma 6:
that the kernel of $\delta_{XM}$ is equal to $\Hom(X,M).$ Namely, if 
$\langle \bdim X,\bDim M\rangle > 0,$ then $\delta_{XM}$ cannot be a monomorphism, thus
the kernel of $\delta_{XM}$ is non-zero. 

	\bigskip\bigskip 
%=========================================================================
{\bf 3. Proof of Theorem.}
	\medskip 
{\bf Lemma 9.} {\it A generic module $M$  has  indecomposable regular
submodules, but no indecomposable preinjective submodules.}
	\medskip
Proof. If $X$ is an indecomposable preinjective module and $N$ is any indecomposable
module with $\Hom(X,N) \neq 0$, then also $N$ is preinjective (see for example [R2]). 
This shows
that a generic module $M$ has no indecomposable preinjective submodule

Let $x = \bDim M$. Then $x$ is a non-zero regular element of $K_0(\Lambda)$, according to 
Lemma 5.
According to Lemma 2, we have $\langle \Phi^{-t}x,x\rangle > 0$ for some $t \ge 0$.
With $x$ also $\Phi^{-t}x$ is a regular element of $K_0(\Lambda)$, thus
according to Lemma 1, there is a regular module $R$ with $\bdim R = \Phi^{-t}x.$
According to Lemma 8, we have $\Hom(R,M) \neq 0$. Let $\phi\:R \to M$
be a non-zero homomorphism and let $X$ be the image of $\phi$. As a factor module of $R$,
the module $X$ is a direct sum of indecomposable modules which are preinjective or
regular. But as we have noted, $M$ has no preinjective submodules, thus $X$ is regular
and of course non-zero. 
	\bigskip
%======================================================
Here is now the proof of Theorem.
Let $M$ be an indecomposable endo-finite module of infinite length.
Let $R(M)$ be the sum of all regular submodules
(since the sum of two regular submodules is again regular,
this actually is the union of all regular submodules).
According to Lemma 9,  we know that $R(M) \neq 0$. Of course, $R(M)$ is a $\Lambda$-submodule:
if $M'$ is a regular submodule and $\phi$ an endomorphism of $M$, then also
$\phi(M)$ is regular. Let us consider the factor module $M/R(M)$. Since
$R(M)$ is both a $\Lambda$-submodule as well as an $E(M)$-submodule, it follows
that $M/R(M)$ is endo-finite.

We claim that any indecomposable submodule of
$M/R(M)$ of finite length is preprojective. Consider a finite length submodule $U$
of $M/R(M)$ which is regular or
preinjective. Then there is a finite length submodule $M'$ of $M$ such that
the canonical map $M' \subset M \to M/R(M)$ maps onto $U$.
Since $R(M)$ is the filtered union of regular modules,
and $M'\cap R(M)$ is a finite length submodule of $R(M)$, there is a regular
submodule $M''$ of $R(M)$ which contains $M'\cap R(M)$. It follows from
$M'\cap R(M) \subseteq M'' \subseteq R(M)$ that $M'\cap M'' = M'\cap R(M)$, thus
$$
 (M'+M'')/M'' \simeq M'/(M'\cap M'') = M'/(M'\cap R(M) \simeq U.
$$
This shows that $M'+M''$ is an extension of the regular module $M''$
by the finite length module $U$ which is regular or preinjective, thus
$M'+M''$ is a direct sum of regular and preinjective modules. Clearly, $M$
has no non-zero preinjective submodules, thus $M'+M'' \subseteq R(M)$. 
But $M' \subseteq R(M)$ implies that the image $U$ of the canonical map
$M' \subset M \to M/R(M)$ is zero. 

Since $M/R(M)$ is endo-finite, it is a direct sum of copies of a finite number
of indecomposable endo-finite modules, say $N_1,\dots,N_t$. 
As we have shown, none of the modules $N_i$ can be regular or preinjective.
Also, if $N_i$ has infinite length, then according to Lemma 9 the module $N_i$ has a 
non-zero regular submodule, but this is impossible. This shows that all the modules
$N_i$ are preprojective. But this implies that  $R(M)$ is a direct summand of $M$,
see for example [R2]. Since $M$ is indecomposable, we conclude that $M = R(M).$
	\medskip
As Krause has pointed out, $R(M)$ is the torsion submodule of $M$ for a torsion pair
and $R(M)$ is a pure submodule of $M$. Thus, if $M$ is indecomposable and endo-finite,
then either $R(M) = 0$ or $R(M) = M.$
	\bigskip\bigskip 
\vfill\eject
%=========================================================================
{\bf 4. Remarks.}
	\medskip 
{\bf (1) A question.} 
It seems to be of interest to determine the set of endo-length vectors 
of the generic modules. One may conjecture that it is the set of all
positive imaginary roots which are not proper multiples of zero roots.
(Here, we call a root $x$ a {\it zero root} provided $\langle x,x\rangle = 0.$)
	\medskip 
%======================================================
{\bf (2) Hereditary artin algebras.} 
The assertions of Theorem and its Corollary are true in the more general setting
of dealing with an arbitrary hereditary artin algebra $\Lambda$, and not just
the path algebras of finite directed quivers. Namely: {\it If $\Lambda$ is a wild connected
hereditary artin algebra, then any generic module is the union of its regular submodules
and the limit closure of the preprojective component does not contain any generic module.}
	
The proof of the general result follows the given one, step by step, only few
alterations are necessary. The actual calculations depend on the decision which kind of 
dimension vectors for finite length modules one wants to use. There are two obvious choices.
Let $\Lambda$ be a $k$-algebra, where $k$ is a commutative artinian ring
such that $\Lambda$ is a finite length $k$-module, and let $M$ be a $\Lambda$-module
of finite length. First of all, one may look at the equivalence class $\bdim M$ of
$M$ in $K_0(\Lambda)$, this means that the coefficient $(\bdim M)_i$ is just the
Jordan-H\"older multiplicity of $S(i)$ in $M$. But one may also take the vector 
$\bdim_k M$ with coefficient $(\bdim_k M)_i$ being the length of $M_i$ as a $k$-module. 
Of course, for the path algebra $\Lambda = k\Delta$ of a directed quiver $\Delta$,
we have $\bdim = \bdim_k$, but in general  $\bdim_k$ is obtained from $\bdim$
by a non-trivial diagonal linear transformation: 
we have to multiply the coefficient $(\bdim M)_i$
by $\dim_k S(i)$. 
	\bigskip
{\bf 5. References.}
	\medskip 
\item{[B]} Baer, D.: Wild hereditary algebras and linear methods. manuscripta math. 55 (1986),
    69-82.
\item{[BGP]} Bernstein, I.~N., Gelfand, I.~M., Ponomarev V.~A.:
  Coxeter functors and Gabriel's theorem. Uspekhi Mat. Nauk 28 (1973), 
  Russian Math. Surveys   28 (1973), 17-32,
\item{[C]}  Crawley-Boevey, W.~W.: Tame algebras and generic modules. 
  Proc. London Math. Soc. 63 (1991), 241-265. 
\item{[K]}  Krause, H.: Abelian length categories of strongly unbounded type.
    arXiv:1301.6665.
\item{[L]} Lenzing, H.: Generic modules over tubular algebras. In: Advances in Algebra
  and Model Theory. Gordon and Breach (1997), 375-385.
%\item{[RR]} Reiten, I., Ringel, C.~M.:
\item{[R1]} Ringel, C.~M.: Representations of K-species and bimodules. 
    J. Algebra 41 (1976), 269-302
\item{[R2]}  Ringel, C.~M.: Infinite dimensional representations of finite dimensional
   hereditary algebras. Symposia Math. 23 (1979), 321-412.
\item{[R3]}  Ringel, C.~M.: The spectrum of a finite dimensional algebra. Proc. Conf.
     Ring  Theory Antwerp 1978. Marcel Dekker Lecture Notes Pure Appl. Math 51 (1979), 535-797. 
\item{[R4]}  Ringel, C.~M.: Tame algebras and integral quadratic forms. Springer Lecture
   Notes in Mathematics 1099 (1984).
	\bigskip

{\rmk
\noindent 
Shanghai Jiao Tong University, Shanghai 200240, P. R. China, and 
  King Abdulaziz University, PO Box 80200,  Jeddah, Saudi Arabia.  E-mail: 
{\ttk ringel\@math.uni-bielefeld.de}
}

\bye